\documentclass[12pt]{article}
\usepackage{amssymb}
\usepackage{amsmath}

\newcommand{\lap}{\mbox{$\bigtriangleup$}}
\newcommand{\grad}{\mbox{$\bigtriangledown$}}
\newcommand{\ra}{{\mbox{$\rightarrow$}}}
\newcommand{\be}{\begin{equation}}
\newcommand{\ee}{\end{equation}}

\newtheorem{mthm}{Theorem}

\newtheorem{rem}{Remark}[section]
\begin{document}

\title{A Hopf type lemma for fractional equations}

\author{Congming Li \thanks{School of Mathematics, Shanghai Jiao Tong University, congmingli@gmail.com, partially supported by NSFC 11571233 and NSF DMS-1405175.}
\quad  Wenxiong Chen \thanks{Partially supported by the Simons Foundation Collaboration Grant for Mathematicians 245486.}}

\date{\today}
\maketitle
\begin{abstract}
In this short article, we state a Hopf type lemma for fractional equations and the outline of its proof. We believe that it will become a powerful tool in applying the method of moving planes on fractional equations to obtain qualitative properties of solutions.
\end{abstract}
\bigskip

\section{Introduction}

The classical Hopf lemma plays a fundamental role in the study of elliptic partial differential equations, it is also a powerful tool in carrying out the method of moving planes to derive symmetry, monotonicity, and non-existence of solutions. Recently, a lot of attention has been turned to the fractional equations due to their broad applications to various branches of sciences. In this paper, we will introduce a fractional version of Hopf type lemma for anti-symmetric functions on half spaces, which can be applied immediately to the second step of the method of moving planes to establish qualitative properties, such as symmetry and monotonicity for solutions of fractional equations.

The fractional Laplacian is a non-local operator defined as
\begin{equation}
(-\Delta)^{\alpha/2} u(x) = C_{n,\alpha} \, \lim_{\epsilon \ra 0} \int_{\mathbb{R}^n\setminus B_{\epsilon}(x)} \frac{u(x)-u(z)}{|x-z|^{n+\alpha}} dz,
\label{op}
\end{equation}
 where $\alpha$ is any real number between $0$ and $2$.

Due to the nonlocal nature of fractional order operators, many traditional methods on local differential operators no longer work. To circumvent this difficulty, Cafarelli and Silvester \cite{CS} introduced the
{\em extension method} which turns a nonlocal problem for the fractional Laplacian into a local one in higher dimensions, hence the classical approaches for elliptic partial differential operators can be applied to the extended equations. This {\em extension method} has been employed by many researchers successfully to obtain interesting results on equations involving the fractional Laplacian. However, so far as we know, besides the fractional Laplacian, there is no {\em extensions methods} that works for other non-local operators, such as the uniformly elliptic non-local operator and fully non-linear non-local operators including the fractional p-Laplacian.
In this paper, we employ a completely different approach and analyze the nonlocal problem directly.
It works not only for the fractional Laplacian, but also for the above mentioned other nonlocal operators as well.

To illustrate how the Hopf lemma can be used in the second step of the method of moving planes, we consider the following simple example
\begin{equation}
(-\lap)^{\alpha/2} u = f(u(x)), \;\; x \in R^n .
\label{ex1}
\end{equation}
Under certain conditions on $f$, we want to show that each positive solution is radially symmetric about some point.

As usual, let
$$T_{\lambda} =\{x \in \mathbb{R}^{n}|\; x_1=\lambda, \mbox{ for some } \lambda\in \mathbb{R}\}$$
be the moving planes,
$$\Sigma_{\lambda} =\{x \in \mathbb{R}^{n} | \, x_1>\lambda\}$$
be the region to the right of the plane, and
$$ x^{\lambda} =(2\lambda-x_1, x_2, ..., x_n)$$
be the reflection of $x$ about the plane $T_{\lambda}$.

Assume that $u$ is a solution of pseudo differential equation (\ref{ex1}). To compare the values of $u(x)$ with $u(x^{\lambda})$, we denote
$$w_{\lambda} (x) = u(x^{\lambda}) - u(x) .$$
Obviously, $w_\lambda$ is anti-symmetric: $w_{\lambda}(x^{\lambda})=-w_{\lambda}(x)$, and it satisfies
$$(-\lap)^{\alpha/2} w_{\lambda}(x) = c_{\lambda}(x) w_{\lambda}(x) , \;\; x \in \Sigma_{\lambda}.$$
The first step is to show that for $\lambda$ sufficiently positive, we have
\begin{equation}
w_{\lambda}(x) \geq 0 , \;\; x \in \Sigma_{\lambda} .
\label{w}
\end{equation}
This provides a starting point to move the plane. Then in the second step, we move the plane to the left as long as inequality (\ref{w}) holds to its limiting position to show that $u$ is symmetric
about the limiting plane. Let
$$\lambda_o = \inf \{ \lambda \mid w_{\mu}(x) \geq 0, \; x \in \Sigma_{\mu}, \, \mu \geq \lambda\}$$
be the lower limit of such $\lambda$ that (\ref{w}) holds.
To show that $u$ is symmetric about the limiting plane $T_{\lambda_o}$, or
\begin{equation}
 w_{\lambda_o}(x) \equiv 0 , \;\; x \in \Sigma_{\lambda_o};
 \label{w0}
 \end{equation}
one usually uses a contradiction argument. Suppose (\ref{w0}) is false, then by a strong maximum principle, we have
$$w_{\lambda_o} (x) > 0 , \;\; \forall \, x \in \Sigma_{\lambda_o} .$$

On the other hand,  by the definition of $\lambda_o$, there exists a sequence $\lambda_k \nearrow \lambda_o$, and $x^k \in \Sigma_{\lambda_k}$, such that
\begin{equation}
w_{\lambda_k}(x^k) = \min_{\Sigma_{\lambda_k}} w_{\lambda_k}  < 0, \;\; \mbox{ and } \; \grad w_{\lambda_k}(x^k) = 0 .
\label{wxk}
\end{equation}

Then under a mild assumption on $f(u)$, we can guarantee that there is a subsequence of $\{x^k\}$ that converges to some point $x^o$, and hence from (\ref{wxk}), we have
$$w_{\lambda_o}(x^o) \leq 0, \mbox{ hence } x^o \in \partial \Sigma_{\lambda_o}; \;\; \mbox{ and } \; \grad w_{\lambda_o}(x^o) = 0 .$$
This will contradicts the following

\begin{mthm} (A Hopf type lemma for anti-symmetric functions)

Assume that $w \in C_{loc}^3(\bar{\Sigma})$,
$$ \overline{\lim_{x \ra \partial \Sigma}} c(x) = o (\frac{1}{[dist(x,\partial \Sigma)]^2}),$$
 and
\begin{equation}
\left\{\begin{array}{ll}
(-\lap)^{\alpha/2}w(x) +c(x)w(x)= 0  &\mbox{ in } \Sigma,\\
w(x) > 0&\mbox{ in }  \Sigma \\
w(x^{\lambda})=-w(x)  &\mbox{ in } \Sigma,
\end{array}
\right.
\label{meq}
\end{equation}
Then
\begin{equation}
\frac{\partial w}{\partial \nu} (x) < 0 \;\;  x \in \partial \Sigma .
\label{partial}
\end{equation}

\end{mthm}
Here, for simplicity of notation,
we denote $\Sigma = \Sigma_{\lambda}$ and $w=w_\lambda$.

\begin{rem}
By a standard interior regularity argument, if $u$ is solution of (\ref{ex1}), even in a week or distributional sense, it is locally smooth provided $f(\cdot)$ is smooth, and consequently, $w_{\lambda}$ is smooth. Hence here the regularity assumption on $w$ can be easily satisfied by
a certain regularity condition on $f$.
\end{rem}
\bigskip

 {\bf Outline of Proof of the Theorem}.

 Without loss of generality, we may assume that $\lambda =0$ and it is suffice to show that $\frac{\partial w}{\partial x_1} (0) > 0$. We argue indirectly. Suppose that $\frac{\partial w}{\partial x_1} (0) = 0$, then it follows from the anti-symmetry of $w$, we also have $\frac{\partial^2w}{\partial x_1^2}(0) = 0$; we will derive a contradiction with the equation.

 The analysis is quite delicate.

 By the definition of the fractional Laplacian, we derive
 \begin{eqnarray}
 & & \frac{1}{C_{n,\alpha}}(-\lap)^{\alpha/2} w(x) \nonumber\\
 &=& \int_{\Sigma} \left\{ \left( \frac{1}{|x-y|^{n+\alpha}} - \frac{1}{|x-y^0|^{n+\alpha}} \right)[ w(x) - w(y) ] +  2 \frac{w(x) d y}{|x-y^0|^{n+\alpha}} \right\} d y \nonumber\\
 &\equiv& \int_{\Sigma} F(x,y) dy + 2 w(x) \int_{\Sigma} \frac{1}{|x-y^0|^{n+\alpha}}  dy \equiv I_1 + I_2.
 \label{1}
 \end{eqnarray}
 Here the first integral is in the Cauchy principal sense.

 We divide $\Sigma$ into several subregions and estimate the above integrals in each region.

 Let
 $$\delta = x_1 = dist(x, T_0), \; x=(x_1, x'),$$
 $$A_{\delta, \epsilon} = \{x \mid 0 \leq x_1 \leq 2 \delta , |x'| < \epsilon \}, \; B_{\epsilon} = \{ x \mid  2 \delta \leq x_1 \leq \epsilon, \, |x'| < \epsilon \} ,$$
$$\Omega_{R, \eta} = \{ x \in B_R(0) \setminus (A_{\delta, \epsilon} \cup B_\epsilon) \mid x_1 > \eta \}, \; D = \{x \mid 1\leq x_1 \leq 2, \, |x'| \leq 1 \}, $$
$$E = \Sigma \setminus (\Omega_{R,\eta}\cup A_{\delta, \epsilon} \cup B_{\epsilon}) .$$
\vspace{4in}

For sufficiently small $\epsilon$ and $\delta$ (with $\delta$ much smaller than $\epsilon$), we estimate $I_1$ on the above sets successively as follows.

1) $$\int_{D} F(x,y) d y \leq -c_1 \delta ,$$
for some $c_1 > 0$ independent of $\epsilon$.

2) $$| \int_{A_{\delta, \epsilon}} F(x,y) d y | \leq c_2 \max\{ \epsilon^{2-\alpha}, \delta^{2-\alpha}, \delta \} \delta. $$

3) $$| \int_{B_{\epsilon}} F(x,y) d y | \leq c_3 \epsilon^{2-\alpha} \delta. $$

Choose $\epsilon$ and $\delta$ small, such that

4) $$| \int_{A_{\delta, \epsilon} \cup B_\epsilon} F(x,y) d y | \leq \frac{c_1}{4} \delta.$$

Fix this $\epsilon$. Let $R$ be sufficiently large and $\eta$ be sufficiently small, such that

5) $$| \int_{E} F(x,y) d y | \leq \frac{c_1}{4} \delta. $$

Fix the above $R$ and $\eta$, let $\delta$ be small, such that $w(x)-w(y) \leq 0, $ for all $y \in \Omega_{R, \eta}$, this is possible because $w(y)$ is bounded away from $0$ in $\Omega_{R, \eta}$ and $w(x) \ra 0$ as $\delta \ra 0$. It follows that
$$\int_{\Omega_{R,\eta}} F(x,y) d y \leq \int_{D} F(x,y) d y \leq - c_1 \delta . $$

Combining the above 5 estimates, we have
\begin{equation}
\int_{\Sigma} F(x,y) d y \leq - \frac{c_1}{2} \delta .
\label{keyE}
\end{equation}

On the other hand, in $I_2$, since $w(x) = O(\delta^3)$, we have
\begin{equation}
w(x) \int_{\Sigma} \frac{1}{|x-y^0|^{n+\alpha}}  dy = O (\delta^{3-\alpha}) \;\; \mbox{ and } \;\; c(x)w(x) \leq o(1) \delta .
\label{keyE1}
\end{equation}

Combining (\ref{1}), (\ref{keyE}), and (\ref{keyE1}), we arrive at
$$(-\lap)^{\alpha/2}w(x) +c(x)w(x) \leq - \frac{c_1 \delta}{4},$$
for $\delta \equiv dist(x, T_0)$ sufficiently small. This
contradicts equation (\ref{meq}). Therefore, we must have $\frac{\partial w}{\partial x_1}(0) = 0.$

\bigskip

{\em Authors' Addresses and E-mails:}
\medskip

Congming Li

School of Mathematical Sciences

Shanghai Jiao Tong University

Shanghai, China, and

Department of Applied Mathematics

University of Colorado,

Boulder CO USA

congmingli@gmail.com
\medskip

Wenxiong Chen

Department of Mathematical Sciences

Yeshiva University

New York, NY, 10033 USA

wchen@yu.edu
\medskip

\end{document}